\documentclass[bsl,bibay2]{asl}

\usepackage[utf8]{inputenc}
\usepackage[french,german,english]{babel}
\usepackage{tikz}

\newcommand{\pind}{\nolinebreak[3]\ensuremath{\mid}}
\newcommand{\makepagination}[1]{\marginpar{\normalfont\scriptsize\sffamily #1}}
\newcommand{\po}[1]{%
\nolinebreak\makepagination{#1}\pind%
}
\newcommand{\ins}[1]{$\langle\!\langle$#1$\rangle\!\rangle$}
\newcommand{\del}[1]{$[\![$#1$]\!]$}
\newcommand{\citep}[1]{[\citeauth{#1}~\citeyear{#1}]}

\makeatletter
\def\citeauth#1{
\nocite{#1}%
\@ifundefined{gu@#1}
   {\cite@nondef{#1}}
   {\csname gu@#1\endcsname}}

\def\citeyear#1{
\nocite{#1}%
\@ifundefined{yr@#1}
   {\cite@nondef{#1}}
   {\csname yr@#1\endcsname\csname yt@#1\endcsname}}
\makeatother

\title[Heinrich Behmann on the decision problem]{Heinrich Behmann's
    1921 lecture on the decision problem and the algebra of logic}


\author{Paolo Mancosu}
\revauthor{Mancosu, Paolo}

\address{Department of Philosophy\\
314 Moses Hall \#2390\\
University of California\\
Berkeley, CA 94720--2390, USA}
\email{mancosu@socrates.berkeley.edu}
\urladdr{http://philosophy.berkeley.edu/mancosu/}

\author{Richard Zach}
\revauthor{Zach, Richard}

\address{Department of Philosophy\\
University of Calgary\\
2500 University Drive N.W.\\
Calgary, AB T2N 1N4, Canada}
\email{rzach@ucalgary.ca}
\urladdr{http://richardzach.org/}

\begin{document}

\begin{abstract}
Heinrich Behmann (1891--1970) obtained his \emph{Habilitation} under
David Hilbert in G\"ottingen in 1921 with a thesis on the decision
problem.  In his thesis, he solved---independently of L\"owenheim and
Skolem's earlier work---the decision problem for monadic second-order
logic in a framework that combined elements of the algebra of logic
and the newer axiomatic approach to logic then being developed in
G\"ottingen.  In a talk given in 1921, he outlined this solution, but
also presented important programmatic remarks on the significance of
the decision problem and of decision procedures more generally. The
text of this talk as well as a partial English translation are
included.
\end{abstract}

\maketitle

\begin{figure}
\includegraphics[width=\textwidth]{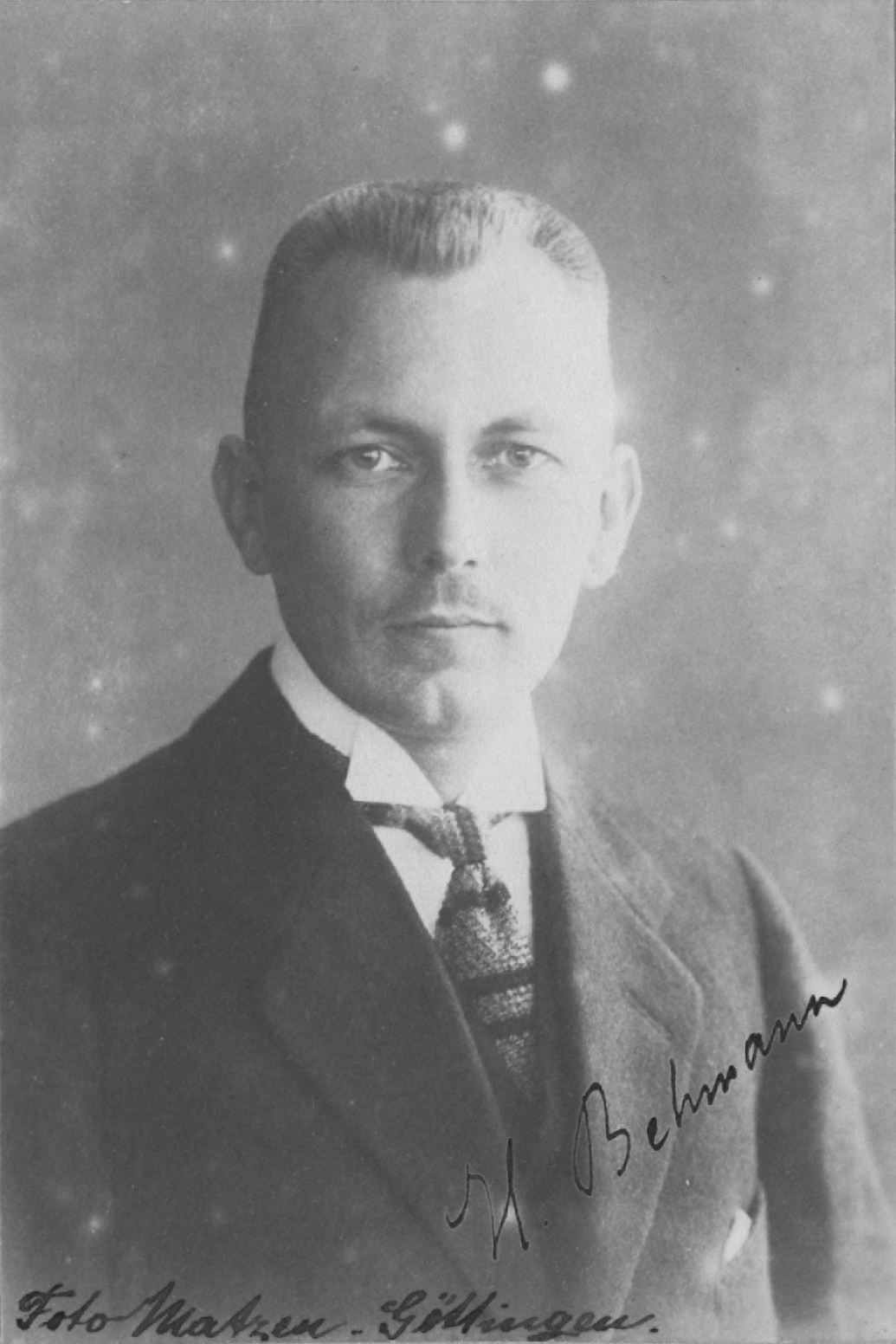}
\end{figure}

\section{Behmann's Career}

Heinrich Behmann was born January 10, 1891, in Bremen. In 1909 he
enrolled at the University of T\"ubingen. There he studied mathematics
and physics for two semesters and then moved to Leipzig, where he
continued his studies for three semesters. In 1911 he moved to
G\"ottingen, at that time the most important center of mathematical
activity in Germany. He volunteered for military duty in World War~I,
was severely wounded in 1915, and returned to G\"{o}ttingen in 1916. In
1918, he obtained his doctorate with a thesis titled \emph{The
  Antinomy of Transfinite Numbers and its Resolution by the Theory of
  Russell and Whitehead} [\emph{\foreignlanguage{german}{Die Antinomie
      der transfiniten Zahl und ihre Aufl"osung durch die Theorie von
      Russell und Whitehead}}] under the supervision of David Hilbert
\citep{Behmann1918}. The thesis is devoted to a reconstruction of the
theory of \emph{Principia Mathematica} motivated by a desire for
removing a number of unclarities as well as by an attempt to put the
project of \emph{Principia} squarely in an empiricist
framework. Recent scholarship
[\citeauth{Mancosu1999}~\citeyear{Mancosu1999},
  \citeyear{Mancosu2003}] has brought to light Behmann's pivotal role
in the understanding and appreciation of \emph{Principia} in
G\"{o}ttingen, thereby providing some essential pieces connecting
Russell's and Whitehead's project with the development of logic and
the foundations of mathematics in Hilbert's school. Behmann then
started to work on the decision problem in logic, for which he coined
the term ``Entscheidungsproblem''. In 1921 Behmann obtained his
\emph{Habilitation} with a thesis titled \emph{Contributions to the
  Algebra of Logic, in particular to the Decision Problem}
     [\emph{\foreignlanguage{german}{Beitr"age zur Algebra der Logik,
           insbesondere zum Entscheidungsproblem}}]. On May 10, 1921,
     Behmann reported his results to the G\"ottingen Mathematical
     Society \citep{Behmann1921}; the text of this lecture is
     reproduced below.  The results were subsequently published in an
     article in \emph{Mathematische Annalen} \citep{Behmann1922}.
     With the \emph{Habilitation} he was conferred the \emph{venia
       legendi}, which allowed him to work and teach as
     \emph{Privatdozent} in the Department of Mathematics in
     G\"{o}ttingen from 1921 to 1925. In 1925 he became professor at
     the University of Halle.

Behmann was not only in close contact with his colleagues in
G\"ottingen (most prominently Hilbert and Bernays, but also with
Wilhelm Ackermann and Moses Sch\"{o}nfinkel) but during the 1920s he
also actively collaborated with a number of people interested in
scientific philosophy such as Carnap and Felix Kaufmann. A visit to
Vienna in 1930 gave rise to an interesting attempt (alas, failed) to
show that classical mathematics is for the most part
intuitionistically justifiable (for the analysis of the full episode,
which also involved G\"odel and Wittgenstein, see
\cite{Mancosu2002a,Mancosu2003}). Behmann is also known for an
interesting solution to the paradoxes of set theory that has been
revived and discussed in the literature on truth (see
\cite{Parsons2003} and \cite{Thiel2002} for the historical
context). On account of his affiliation with the NSDAP he was
dismissed from his position in Halle in 1945. Fearing deportation, he
fled Soviet-occupied Saxony to his home town Bremen in 1946. Due in
part to his wartime injury, he was unable to obtain academic
employment. He died in Bremen on February 3, 1970. Behmann's extensive
\emph{Nachla\ss} was for many years in Aachen and later Erlangen, and
is now housed at the Staatsbibliothek zu Berlin. See
\cite{HaasStemmler1981} for additional detail on Behmann's life and
work and an overview of Behmann's papers.

\section{The Decision Problem: from the Algebra of Logic to 
  \textit{Principia}}

The reference to the ``decision problem'' in the title of Behmann's
lecture may well be the first time the expression made its public
appearance. Behmann claimed to have coined the word on several
occasions. For instance, he wrote to Russell in 1922:
\begin{quote}
  It was what I call the Problem of Decision, formulated in the said
  paragraph [\P1 of \cite{Behmann1922}] that induced me to study the
  logical work of Schr\"{o}der. And I soon recognized that, in order
  to solve my particular problem, it was necessary first to settle the
  main problem of Schr\"oder's Calculus of Regions, his so-called
  Problem of Elimination. [\dots] Indeed, the chief merit of the said
  problem [the Decision Problem] is, I daresay, due to the fact that
  it is a problem of fundamental importance on its own account, and,
  unlike the applications of earlier Algebra of Logic, not at all
  imagined for the purpose of symbolic treatment, whereas, on the
  other hand, the only means of any account for its solution are
  exactly those of symbolic logic. (Behmann to Russell, August 8,
  1922. Russell Archive, McMaster University)
\end{quote}

In a letter to Scholz (December 27, 1927, Behmann Archive,
Staatsbibliothek zu Berlin), Behmann claims to have been the first to
``explicitly'' [\emph{\foreignlanguage{german}{ausdr"ucklich}}] pose
the general decision problem. The \emph{general} decision problem is
the extension of the decision problem from propositional logic to full
logic. In the same letter, he grants that related problems had emerged
in the algebra of logic tradition but that the full problem had not
been stated before him. He also stresses that one should not confuse
this ``Entscheidungsproblem'' with the ``Entscheidungsproblem''
formulated, in connection to Grelling's paradox, by
\cite{Hessenberg1906}.  Behmann also points to the connection between
the decision problem and the ``elimination problem'' central to the
algebra of logic. He says that the elimination problem must be
considered to be independent of the general decision problem while it
is at the same time a special case of it. Behmann mentions Peirce,
Schr\"oder, and L\"owenheim as among the most important contributors
to work on the elimination problem. (Here it is important to point out
that Behmann achieved his important results on the decidability of the
monadic second-order predicate calculus independently of L\"owenheim,
whose work he discovered after he had found the result.)  In the same
letter to Scholz, Behmann mentioned that the he had learned the
solution of the decision problem for propositional logic using normal
forms directly from Hilbert.\footnote{Already \cite{Hilbert1905a}
  considered conjunctive and disjunctive normal forms and decidability
  of propositional logic.  The first explicit discussion of the
  decision problem for propositional logic was given in
  \cite{Bernays1918}, but the essential idea is found in
  \cite{Hilbert1917}; see \cite{Zach1999} for discussion).}

Moving now to the lecture, let us begin by pointing out that in his
[\citeyear{Behmann1922}], Behmann made explicit what in his lecture he
only presupposes, namely that ``As specifically Peano and Russell have
seen, every mathematical proposition can fundamentally be thought of as
a purely logical state of
affairs.''\footnote{``\foreignlanguage{german}{Wie namentlich Peano
    und Russell erkannt haben, l"a"st sich in der Tat jeder
    mathematische Satz letzlich als ein rein logischer Sachverhalt
    auffassen}'' (\cite{Behmann1922}, p.~166, note~1).}  Since, thanks
in particular to the work of Whitehead and Russell, all of logic can
be axiomatized, it follows from the expressibility of every
mathematical statement as a logical one that every mathematical
question could be decided if it were possible to decide for every
logical sentence whether it could be derived from the axioms. And this
leads to the claim, made in Behmann's lecture in 1921 and in
\cite[p.~166]{Behmann1922}, that a solution of this problem would
transform mathematics into ``one enormous
triviality.''\footnote{\cite{Bernays1918} also suggested that the
  decision procedure for propositional logic reduces it to triviality.
  The idea that a decision procedure for all of logic would trivialize
  mathematics led von \cite{Neumann1927}, on the other hand, to doubt
  that such a procedure can be found (cf.~\citep{Gandy1995},
  \S5.2). \cite{Weyl1927} similarly writes, ``Completeness in this
  sense would only be ensured by the establishment of such procedural
  rules of proof as would lead demonstrably to a solution for every
  pertinent problem. Mathematics would thereby be trivialized. But
  such a philosopher's stone has not been discovered and never will be
  discovered'' (p.~20; p.~24 of the translation).  Behmann also
  compares a decision method to the ``philosopher's stone'' in the
  lecture.}  To repeat, while the connection is left implicit in the
1921 lecture, it is also obviously presupposed in the way the lecture
moves from a reflection on mathematics to the full decision
problem. In the lecture, Behmann continues by asserting that the fact
that logic (and hence mathematics) can be axiomatized is not enough to
turn mathematics into a triviality, for the rules of derivation tell
us ``only what one may do, and not what one should do.''  Behmann then
describes the decision problem as the more specific problem of finding
a deterministic, computational procedure to decide any mathematical
claim:
\begin{quotation}
[We require] not only the individual operations but also the
\emph{path of calculation} as a whole should be specified by rules, in
other words, \emph{an elimination of thinking in favor of mechanical
  calculation.} If a logical or mathematical assertion is given, the
required procedure should give complete instructions for determining
whether the assertion is \emph{correct or false by a deterministic}
[\emph{\foreignlanguage{german}{zwangsl"aufig}}] \emph{calculation
  after finitely many steps}.  The problem thus formulated I want to
call the \emph{general decision problem}.
\end{quotation}
While Behmann's programmatic remarks are noteworthy for their analysis
of the significance of the decision problem, and for his emphasis on
deterministic computation, neither they nor the methods for solving
the decision problems developed by him yet amount to a coherent
analysis of the notion of deterministic computational procedure.

In order to make headway on the decision problem for logic, Behmann
suggests that the axiomatic approach is unsuitable, and that one
instead has to adopt the perspective of the algebra of logic. As
mentioned above, Behmann saw the decision problem as connected to the
elimination problem in the algebra of logic tradition, and so he found
it necessary to work through it.  We cannot here give a full
explanation of the elimination problem in the algebra of logic, which
is itself only a sub-problem of what is usually called the ``solution
problem'' [\emph{\foreignlanguage{german}{Aufl"osungsproblem}} in
  Schr\"oder]. The latter problem was already stated by
Boole. Basically, the idea is that from a set of logical equations
(usually containing unknowns) one would like to find (some or all)
``solutions'' for such equations, namely values of the variables that
make the equations true. In the hands of Schr\"{o}der, this process
leads to the study of the ``resultants'' $r_P(x_1, \dots, x_n)=0$ of
an equation which eliminates a variable~$y$ from an equivalent
polynomial inequality of the form $P(x_1, \dots, x_n, y) \neq 0$. The
``elimination'' problem is the determination of such ``resultants''
(for a detailed analysis and further references see \cite{Bondoni2007,
  Bondoni2009}). Rather than attempting a presentation of the context
of the two problems within the algebra of logic tradition, let us
mention how they emerge within the second-order predicate calculus
with binary predicates (the restriction to unary predicates or the
generalization to $n$-ary ($n > 2$) predicates is immediate). In the
decision problem one only investigates the validity, or lack thereof,
of any well-formed formula~$A$ of the calculus. In the elimination
problem, one also asks for a well-formed formula~$B$ that is logically
equivalent to~$A$ but simpler (thereby simplifying the determination
of its validity or lack thereof). For instance, given the formula
\[
\exists X \forall x\forall y((Xxy \lor Axy) \land (\lnot Xxy \lor Bxy))
\]
we might want a first-order formula equivalent to it. It is not too
hard to check that
\[
\forall x\forall y(A(x,y) \lor B(x,y))
\]
is one such formula. This is the ``resultant'' of the elimination
problem (see \cite{Ackermann1934}). According to Ackermann, the
decision problem is actually a special case of the elimination
problem. Behmann's result provides a solution to the elimination
problem for formulas containing monadic predicate variables and this
immediately yields the solution to the decision problem for the same
class of formulas. Neither the elimination nor the decision problem
for variable binary predicates admits of a general solution. We now
describe Behmann's result.

\section{Behmann's Result}

Behmann's result was first announced in his May 1921 lecture and
published the next year in \emph{Mathematische Annalen}
\citep{Behmann1922}. Behmann also reported it in September 1923 at
the annual meeting of the German Mathematical Society in Marburg
\citep{Behmann1923}. The result gives a method for deciding, given a
formula of monadic second-order logic, whether it is valid
[\emph{\foreignlanguage{german}{allgemeing\"ultig}}]. Behmann
considers three classes of formulas: domain $A$ is the class of
first-order formulas with monadic predicate symbols (but without
identity), domain~$B$ is the class of second-order formulas with
monadic predicate variables (but without identity), and $B^*$ is $B$
plus the identity predicate. Behmann's aim is to solve the elimination
problem for~$B^*$, i.e., to show that any monadic second-order formula
is equivalent to a first-order formula.  Behmann did not have a
precise semantics for his languages available, so the solution to the
decision problem, for him, involved finding a procedure to determine
whether a given formula was true under its intuitive interpretation
for every possible domain.  Since formulas with predicate constants do
not express determinate propositions independently of an
interpretation of the predicate symbols, Behmann's notion of validity
properly applies only to second-order sentences containing no
predicate constants. For instance, Behmann considers the decision
problem to apply not to formulas like $\forall y(Ay \lor \lnot Ay)$,
but only to sentences such as $\forall X\forall y(Xy \lor \lnot Xy)$.
Solving the elimination problem for such sentences, then, would show
that they are each equivalent to a first-order sentence involving no
extralogical symbols, only identity, and hence expresses a condition
on the size of the domain.  The sentence is valid if and only if the
condition is compatible with any size of the domain whatsoever.

In keeping with the tradition of using normal forms in the
investigation of questions like the decision problem, Behmann
considers which special forms formulas can be brought into.  If all
predicate symbols are monadic, then every formula is equivalent to a
formula which is a propositional combination of formulas of the form
$\forall x(A_1x \lor \dots \lor A_nx)$ and $\exists x(A_1x \land \dots
\land A_nx)$, where $A_ix$ is atomic or the negation of an atomic
formula.  

Second-order quantifiers $\forall X$ and $\exists X$ distribute over
$\land$ and $\lor$, respectively, just like first-order quantifiers
do.  Behmann is thus able to reduce the problem to that of removing
the second-order quantifier from a formula of the form $\exists X\,
\phi$ where $\phi$ is a conjunction of formulas of the forms $\forall
y(Fy \lor Xy)$, $\forall y(Fy \lor \lnot Xy)$, $\exists y(Fy \land
Xy)$, and $\exists y(Fy \land \lnot Xy)$. Each of these formulas can
be translated into Behmann's own version of class notation (this
notation expresses certain relationships between the classes denoted
by $F$ and $X$, e.g., $\exists y(Fy \land Xy)$ expresses that $F \cap
X \neq \emptyset$).  Behmann's notation is more flexible than the
class-calculus notation of Schr\"oder, who had been unable to solve
the elimination problem even for the monadic case.  However, the
first-order equivalents obtained by Behmann for formulas in ``main
elimination form'' in general will contain identity.  In the last part
of his paper, then, Behmann generalizes the procedure to formulas of
domain $B^*$. His generalized class notation then contains as basic
expressions formulas that state that the intersection of several
classes contains at least $n$ objects, and its dual, that the union of
several classes contains all \emph{but} $n$ objects.  The elimination
problem is solved by showing that a second-order quantifier can be
removed from a conjunction (in the case of $\exists X$) or disjunction
(in the case of $\forall X$) of such expressions. If the original
second-order formula contains no predicate or individual constants,
the resulting formula is a propositional combination of statements of
the form ``there are at least $n$ objects'' and ``there are at most
$n$ objects''.

\section{Further Developments}

For anyone familiar with the standard textbook proof of the
decidability of monadic logic, Behmann's proof seems terribly
unwieldy---and surprisingly so, since the usual method is so simple.
The usual method consists in showing that if the sentence has a
countermodel, it has a finite countermodel of a size depending on the
number of predicate symbols in the sentence (specifically, $2^k$ for
$k$ the number of predicate symbols).  Since one can check all such
potential countermodels, the question is decidable. This method of
solving the decision problem for monadic sentences was first outlined
by \cite{Hilbert1922b} in a lecture course in the Winter semester of
1922/23.  Parts of the relevant sections of the lecture notes were
later incorporated almost verbatim into the textbook by
\cite{Hilbert1928b}.

In the same semester, the Russian mathematician Moses Sch\"onfinkel
outlined a decidability proof for validity of sentences with a single
binary predicate symbol in prenex form with the quantifier prefix
$\exists x\forall y$. The manuscript containing the proof remains
unpublished, but Bernays published a paper containing the result,
extended to sentences with an arbitrary number of binary predicate
symbols, in 1927.  The result was generalized a year later by
\cite{Ackermann1928}, who showed that validity of formulas with
quantifier prefix $\exists^*\forall\exists^*$ is decidable. The paper
by \cite{Bernays1928} is notable, then, not for the proof of
decidability of the special case of the Ackermann class, but for three
other contributions.  The first is the observation that instead of
deciding validity of formulas, it is easier to focus on
\emph{satisfiability}---a concept that Bernays likely adopted from
Skolem and L\"owenheim.  Note that no definition of model was
available yet, and neither validity nor satisfiability were defined
precisely.  The second was the first proof in print of the
decidability of satisfiability of monadic sentences by giving a bound
on the size of models.  The same approach is used in most subsequent
papers on the decision problem. Indeed, already Ackermann's paper the
following year established the result mentioned by proving that the
class of prenex sentences of the form $\forall^*\exists^* A$ is
``finitely controllable.''  Bernays's third contribution---and the
contribution for which the paper is mainly known---is the proof of the
decidability of satisfiability of sentences in the so-called
Bernays-Sch\"onfinkel class (prenex sentences of the form
$\exists^*\forall^* A$).

Hilbert mentioned the decision problem as one of the main open
problems in the foundations of mathematics in his address to the 1928
International Congress of Mathematicians in Bologna
[\citeauth{Hilbert1928}~\citeyear{Hilbert1928},
  \citeyear{Hilbert1928a}]; it is also mentioned prominently in the
textbook by \cite{Hilbert1928b}.  Subsequently, a number of other
logicians not directly connected to Hilbert's group in G\"ottingen
contributed partial solutions to the decision problem.  These include
Ramsey, Herbrand, K\'almar, and G\"odel. A good presentation of the
solvable cases of the decision problem up to the 1950s is given by
\cite{Ackermann1954}. More comprehensive treatments can be found in
\cite{Dreben1979} and \cite{Borger1997}.

The general case of the decision problem for logic was shown to be
unsolvable in 1935--36, independently by \cite{Church1936} and
\cite{Turing1937}.  This result relied not only on the work of
\cite{Godel1931}, but it required in particular mathematically precise
and general definitions of computational procedures, in particular,
decision procedures.  Church's notion of lambda-definable functions
and Turing machines provided some of the first such definitions. They
were anticipated by related work carried out by Emil Post in the
1920s, who, however, never published it.\footnote{\cite{Post1921a}
  does mention a decision procedure for a ``deductive system involving
  primitive propositions of but one argument,'' but this system was
  propositional only.  See \cite{Urquhart2009} for discussion of
  \cite{Post1965}, and \cite{Gandy1995} for a historical survey of
  the development of notions of computability.}  Behmann's
contribution to the decision problem, like other algorithmic methods
for solving mathematical problems, was in the first instance an
example of such a computational decision procedure. The significance
of Behmann's lecture lies mainly in the clear articulation of the very
general nature of the problem to be decided---essentially, every
mathematical question---and his incisive reflections in the
programmatic part of his lecture on the nature of computational
processes.

\section*{Acknowledgments}\normalfont

We would like to thank the Staatsbibliothek zu Berlin for permission
to publish Behmann's lecture, Prof.~Christian Thiel for originally
providing access to Behmann's papers, and Johannes Hafner for his help
in typesetting the German text.  Behmann's portrait is held by the
Staats- und Universit\"atsbibliothek G\"ottingen, Cod.~Ms.\ Hilbert
754, Bl.~22, Nr.~114, and reproduced here by their kind
permission. Suggestions and comments by Patricia Blanchette, Christian
Thiel, and above all an anonymous referee, have improved the
presentation. Richard Zach would like to acknowledge the support of
the Social Sciences and Humanities Research Council of Canada.

\bibliographystyle{asl} \bibliography{history-of-logic}

\section*{Decision problem and algebra of logic}\normalfont

The translation of the first seven pages of Behmann's lecture is by
Richard Zach.  For readability, underlining in Behmann's manuscript
has been incorporated (as italics) only sparingly, and editorial
emendations other than pagination has been omitted in the translation.
The complete text in the original German follows below.\bigskip

\noindent\textit{Decision problem and algebra of logic}
\hfill \mbox{May 10, 1921}\bigskip

\noindent A couplet, which goes something like this:

\begin{quote}
Sunt mathematici veri natique poeta;\\
sed, quae finxerant, hos probare decet.
\end{quote}
is due, as far as I know, to Kronecker. In translation: 
\begin{quote}
We mathematicians are true poets with a calling;\\
But we still must prove our poetry!
\end{quote}

And we will have to admit that there is something true in these words.
The work of the mathematician can indeed be compared with that of the
poet, or generally with that of a creative artist. They, too, need
\emph{imagination}, in order to obtain a sufficiently large range of
possibilities and to focus on certain conceivable connections between
facts or between thoughts for carrying out their investigations; they
need \emph{trained sensitivity} for identifying the correct one among
a dizzying array of possibilities, and for making lucky guesses of so
far unsuspected states of affairs; they need \emph{subtlety} to
present their results in their suitable form.  As Boltzmann once said,
``among all artists, the mathematician is closest to the demiurge.''

At the same time, the mathematician is not just a creative artist,
but---and this tears us at once from our aesthetic speculation---they
are subject to the additional requirement,\,\po{2}\,to which Kroneker
points in his couplet, that their creations---if we can call it
that---must withstand the sharpest test of critical reason.  They
ought not just to \emph{assert} theorems, even if they are beautiful
and correct, but they must also \emph{prove} them, and with it add
them to the store of secured knowledge.  No other science sets up this
requirement with the same inexorable rigor, and anyone who understands
the nature of mathematics will not attempt to haggle to weaken it, but
instead to tighten it as much as possible.

There is now no doubt that in practice this very same requirement is
an extremely significant burden for the mathematician.  The reason for
this is roughly the following:

While in the beginning mathematics---and this is not just true for
mathematics as a whole but for every single one of its disciplines,
regardless of how early or late they may have come to be
developed---always first focuses on \emph{specific single problems}
and attempts to solve them with whatever methods it already has at its
disposal, it nevertheless, in the course of its development,
eventually develops \emph{general methods} for \emph{wide classes} of
problems (even if these methods do not solve every problem of a
certain kind). The problems to which these methods apply then do not
require any \emph{inventive insight} but only the calculatory tools
provided.  Once this is accomplished, the path to a solution no longer
goes by \emph{searching} and \emph{cutting a trail} in every direction
of the compass, but it is now \emph{completely given},\,\po{3}\,fully
developed, and furnished with clear markers; now the path, as I would
like to put it, \emph{inevitably} leads to the destination.  The
\emph{problem} has thus turned into a mere \emph{exercise in
  computation}.

Unfortunately the kind of progress I have just characterized is
nowhere to be found when the problem is that of \emph{proving
  mathematical propositions}.  The attitude long left behind in the
case of mathematical problems in general still applies in the case of
proving mathematical theorems: it indeed requires a \emph{new,
  creative thinking for each case}, a lucky insight which of course
can never be forced even with the most enormous of efforts.  The task
of \emph{verifying the validity of a given mathematical claim} at
first presents itself as being of a very different nature than that of
solving a cubic equation or of carrying out the Gaussian interpolation
algorithm. For the latter we know \emph{that} we can always find a
solution as well as \emph{how} to find it, for the former neither one
nor the other.  For the latter, we have fully developed,
perspicuous and deterministic procedures, for the former the important
aspect requires instinct, a serendipitous intuition. We are thus
unfortunately still far from being able to delegate the work of
proving mathematical theorems to mathematical laborers, in the way we
are regarding numerical calculations.

At first glance the development of a general, deterministic procedure
for proving mathematical propositions may seem like a \emph{hopeless,
  indeed ludicrous endeavor},\,\po{4}\,perhaps comparable to the
medieval quest for the philosopher's stone. And indeed the achievement
of this goal would be no less significant than the wonders believed to
be inherent in that mythical substance for not just mathematics but
for knowledge as a whole. If we look more closely, the entirety of
mathematics would indeed be transformed into one \emph{enormous
  triviality}; we would be able to solve every problem posed to
us---perhaps not in the intended sense, which may be impossible---,
but at least in the sense in which it can be solved at all, as set out
by Hilbert in his Paris address ``Mathematical Problems.''

From a purely mathematical perspective, however, there does indeed
appear to be no possibility for attacking this problem in its enormous
variety, say, by dividing mathematical propositions into classes
according to their complexity and by advancing from the simple to the
more difficult.

Nevertheless we have, or so I would like to claim, attained a level of
knowledge in an area which is not directly mathematical, but indeed
\emph{much more fundamental}, that we \emph{not only have the ability}
but \emph{every right} to formulate this problem explicitly and
rigorously and to begin its investigation. Indeed we must
\emph{descend to a lower level}, since we are really no longer working
on a properly mathematical problem, but a ``supramathematical''
problem. This lower level is \emph{symbolic logic}, which has been
anticipated by Leibniz, first taken on by Boole,\,\po{5}\,and
developed extensively by Schr\"oder, Frege, Peano, and Russell.

In fact, at its present state of development, it already provides us
with the ability to formulate the problem in a relatively simple and
transparent form, by allowing us to express all propositions in
questions in a uniform way using just a few symbols. Specifically, we
can express \emph{every logical or mathematical proposition with the
  following few words alone:}

\bigskip
\begin{raggedright}
\begin{tabular}{@{}cccccc@{}}
& & & & \multicolumn{2}{c}{$\overbrace{\hspace{1.5in}}^{\text{Begriff}}$}\\
\emph{(proposition)} & \emph{not} & \emph{or} & \emph{object} & \emph{property} & \emph{relation} \\
$p$ & $\overline{p}$ & $p \, q$ & $a$ & $f^{x}$, $F^{\phi}$ & $f^{xy}$ 
\end{tabular}\hfil\par
\hfill\begin{tabular}{@{}ccc@{}}
\emph{satisfies} & \emph{all} & \emph{identical} \\
$f_{a}$, $f_{ab}$ & $xf_{x}$ & $a=b$
\end{tabular}
\end{raggedright}
\bigskip

Of course this does not mean that in practice one always has to carry
out the reduction this far; instead, one will introduce new signs, as
far as is desirable for the particular application, for other
concepts, such as \emph{and, there are, the class, the numbers, the
  operation, etc.}

With this insight we already have made progress.  We now only have to
consider \emph{expressions formed according to rules from these
  symbols} and thus have cast off the bonds of natural language.

As is well known, symbolic logic can be \emph{axiomatized}, i.e.,
reduced to a system of relatively few basic formulas and rules, so
that proving theorems now appears to be a mere computational procedure.
One only has to write down new formulas to formulas already given,
where the rules already specify what may be written down in every
case.  Proving then has taken on the \emph{nature of a game}.  It is
similar to,\,\po{6}\,say, chess, where by moving one's own pieces,
perhaps while removing one of the opponent's pieces, one transforms a
given position into a new one, where only the movement of one's own
piece and the removal of the opponent's must be allowed by the rules
of the game.

But this very comparison brings out quite starkly that the view
of symbolic logic just detailed cannot suffice at all for our problem.
For it shows us, like the rules of chess, \emph{only what one may do},
and \emph{not what one should do}.  The latter remains---in the one as
in the other case---a question of inventive \emph{thinking}, of lucky
\emph{combination}.  We, however, require a lot more: not only the
individual operations but also the \emph{path of calculation} as a
whole should be specified by rules, in other words, \emph{an
  elimination of thinking in favor of mechanical calculation.}  If a
logical or mathematical statement is given, the required procedure
should give complete instructions for determining whether the
statement is \emph{correct or false by a deterministic calculation
  after finitely many steps}.  The problem thus formulated I want to
call the \emph{general decision problem}.

It is essential to the character of this problem that as method of
proof only \emph{entirely mechanical calculation} according to given
instructions, without any activity of thinking in the narrower sense,
is allowed. One might, if one wanted to, speak of \emph{mechanical} or
\emph{machine-like thinking}.  (Perhaps one can one day even let it
be carried out by a machine.)~\po{7}

It is thus not merely symbolic logic as such that is required, but a
\emph{specific direction within symbolic logic}, which not only
investigates how one \emph{may} compute, but how one \emph{should}
compute in order to attain a given goal.  For a particular historical
reason I would like to call this direction that of the \emph{algebra
  of logic}, specifically because all important preparatory work in
this area are subsumed under this name of algebra of logic (Boole,
De~Morgan, Peirce, Schr\"oder). This subsumption is perhaps only
accidental, since none of those just mentioned have seriously
considered an application to the decision problem.  But also in
regards to its substance the name seems to me to be somewhat
appropriate, since this narrower domain will have precisely the
development of \emph{algorithms}, i.e., deterministic procedures, such
as those in algebra for finding solutions to equations, as its aim. I
consider it a fortuitous coincidence that here we are afforded an
opportunity to summarize these, in their superficial form relatively
dubious and obscure---if not to say: boring and
tedious---investigations under a \emph{new, unified and valuable
  perspective}, and with it to save a large amount of deep and arduous
intellectual labor from the inglorious fate of slowly mouldering away
in libraries.  For, as far as I can tell, general interest in the
investigations of that older algebra of logic has almost completely
ceased even among proponents of symbolic logic.~\po{8}

I do not have to emphasize that despite all this we are today still
\emph{extremely far away} from a solution to the decision problem
formulated above.  At least I have succeeded in completely carrying
through this problem for a certain \emph{modest domain of
  propositions}, which nevertheless already contains quite a
significant variety of propositions.  It is this work, on which I
would like to speak to you in particular.

\section*{Entscheidungsproblem und Algebra der Logik}\normalfont

The following manuscript by Heinrich Behmann, entitled
``\foreignlanguage{german}{Entscheidungsproblem und Algebra der
  Logik},'' is the text to a lecture Behmann gave in the G\"ottingen
Mathematical Society on May 10, 1921. The date and venue of the
lecture are attested in the list of events at the \emph{Mathematische
  Gesellschaft in G\"ottingen} in {\bfseries\itshape Jahresbericht der
  Deutschen Mathematiker Vereinigung}~30 (1921), 2.~Abteilung, p.~47.
The manuscript is part of the Behmann papers held at the
Staatsbibliothek zu Berlin, Preu\ss{}ischer Kulturbesitz,
Handschriftenabteilung under call number Nachl.~355 (Behmann), K.~9
Einh.~37, and is reproduced here by their kind permission.

Underlining in red pencil in the original is rendered as italics,
except for formulas in the running text. Words or passages crossed out
by Behmann are enclosed in double square brackets,
\del{\dots}. Authorial insertions are enclosed in double angle
brackets \ins{\dots}. Words or phrases enclosed in single angle
brackets $\langle\dots\rangle$ and single square brackets $[\dots]$ by
Behmann himself are displayed in the same way in the
transcription. Page breaks are indicated by\,\po{1}\,in the running
text and the page number in the margin.\bigskip

\selectlanguage{german}
\noindent\emph{Entscheidungsproblem und Algebra der Logik}
\hfill \mbox{10. Mai 1921.}\bigskip

\noindent Von \emph{Kronecker} stammt, soviel ich wei"s, ein \emph{Distichon},
das ungef"ahr so lautet:

\begin{quote}
Sunt mathematici veri natique poeta;\\
sed, quae finxerant, hos probare decet.
\end{quote}

In "Ubersetzung:

\begin{quote}
Wir Mathematiker auch sind echte, berufene Dichter;\\
uns liegt noch der Beweis f"ur das Gedichtete ob!
\end{quote}

Und man wird zugeben m"ussen, da"s \emph{etwas Wahres} in diesen
Worten liegt. Die \emph{Arbeit des Mathematikers} l"a"st sich in der
Tat bis zu einem gewissen Grade \ins{wohl} mit der des \emph{Dichters}
oder allgemeiner "uberhaupt des \emph{schaffenden K"unstlers}
vergleichen. Auch er bedarf der \emph{Phantasie}, um f"ur die
Durchf"uhrung seiner Untersuchung von vornherein einen hinreichenden
Kreis von M"oglichkeiten zu gewinnen und gewisse denkbare Tatsachen-
oder Gedankenzusammenh"ange ins Auge zu fassen, des \emph{geschulten
  Gef"uhls}, um aus der mitunter verwirrenden Vielzahl von
M"oglichkeiten gerade die richtige herauszufinden, bisher nicht
vermutete Sachverhalte gl"ucklich zu erraten, des \emph{feinen
  Taktes}, um seine Ergebnisse in der ihnen angemessensten Form
darzustellen. Wie \emph{Boltzmann} einmal gesagt hat, "`\emph{kommt
  der Mathematiker unter allen K"unstlern dem Weltensch"opfer am
  n"achsten}."'

Indessen: der Mathematiker ist nun nicht blo"s schaffender K"unstler,
sondern~--- und dies rei"st uns mit einem Schlage aus unserer
"asthetischen Betrachtung heraus~--- f"ur ihn besteht noch die
\emph{zus"atzliche Forderung},\,\po{2}\,auf die \emph{Kronecker} in
seinem Distichon hinweist, da"s seine "`Sch"opfung"'~--- wenn wir so
sagen d"urfen~--- auch die \emph{sch"arfste Probe des kritischen
  Verstandes} aushalten mu"s. Er soll nicht blo"s S"atze
\emph{aufstellen}, und seien sie noch so sch"on und richtig, sondern
er soll sie auch \emph{beweisen} und damit erst dem bisherigen
gesicherten Erkenntnisbesitz einf"ugen. In keiner anderen Wissenschaft
wird diese Forderung mit der \emph{gleichen unerbittlichen Strenge}
gestellt, und wer das Wesen der Mathematik recht versteht, sucht von
dieser Forderung \emph{nicht etwas abzuhandeln}, sondern sie im
Gegenteil so weit wie irgend tunlich noch zu \emph{versch"arfen}.

Es unterliegt nun aber keinem Zweifel, da"s eben diese Forderung in
der Praxis eine ganz \del{\emph{au"serordentliche}}
\ins{\emph{bedeutende}} \emph{Belastung} f"ur
den Mathematiker bedeutet. Und zwar ist der \emph{Grund hierf"ur}
wesentlich der folgende:

W"ahrend anf"anglich die \emph{Mathematik}~--- und zwar gilt dies
nicht nur f"ur die Mathematik im ganzen, sondern f"ur \emph{jedes
  einzelne ihrer Teilgebiete}, so fr"uh oder so sp"at dieses auch zur
tats"achlichen Entwicklung gekommen sein mag~--- immer zun"achst
\emph{bestimmte einzelne Probleme} ins Auge gefa"st und diese nun auf
irgend eine Weise mit den ihr bereits zu Gebote stehenden Mitteln zu
bew"altigen sucht, gelangt sie dessenungeachtet im Verlaufe ihrer
Entwicklung mehr und mehr dahin, wenn auch nicht f"ur alle Probleme
von einer irgendwie bestimmten Natur, so doch f"ur \emph{ausgedehnte
  Klassen} solcher Probleme \emph{allgemeine Verfahren} zu entwickeln,
soda"s es f"ur die von diesen Verfahren beherrschten Probleme dann
\emph{keines erfinderischen Scharfsinns}, sondern nur noch der
Beherrschung der vorausgesetzten rechnerischen Hilfsmittel
bedarf. Damit ist der Weg zur L"osung \emph{nicht mehr} ein nach
irgend einer Himmelsrichtung \emph{zu suchender} und neu \emph{zu
  bahnender}, sondern er ist nunmehr ein \emph{fertig
  vorhandener},\,\po{3}\,voll ausgebauter und mit deutlichen
Wegweisern versehener, der, wie ich sagen m"ochte, \emph{zwangl"aufig}
auf das Ziel hinf"uhrt. Das \emph{Problem} ist damit gewisserma"sen
zur blo"sen \emph{Rechenaufgabe} geworden.

Von einer Entwicklung, wie ich sie eben kennzeichnete, ist nun in dem
Falle, da"s es sich um das \emph{Beweisen mathematischer Aussagen}
handelt, leider \emph{nicht im entferntesten die Rede}. Was bei den
mathematischen Problemen im allgemeinen l"angst "uberwundener
Standpunkt ist: beim Beweisen mathematischer S"atze bedarf es in der
Tat \emph{von Fall zu Fall des neuen erfinderischen Gedankens}, einer
gl"ucklichen Eingebung, die sich bekanntlich auch durch die
gewaltigste aufgewendete M"uhe \emph{niemals erzwingen} l"a"st. Die
Aufgabe, \emph{eine gegebene mathematische Behauptung auf ihre
  G"ultigkeit zu pr"ufen} stellt sich ihrer Natur nach zun"achst
jedenfalls als eine ganz andere dar, als etwa diejenige, \emph{eine
  kubische Gleichung aufzul"osen} oder \emph{eine Gau"ssche
  Ausgleichungsrechnung durchzuf"uhren}. Bei der \emph{zweiten} wissen
wir, \emph{da"s} wir die L"osung unbedingt finden k"onnen und
\emph{wie} wir sie finden k"onnen, bei der \emph{ersten} weder das
eine noch das andere. Bei der \emph{zweiten} sind wir im Besitz
vollst"andig ausgebauter, "ubersichtlicher und \emph{zwangl"aufiger
  Verfahren}, bei der \emph{ersten} bleibt das Wesentliche dem
\emph{Gef"uhl}, dem gl"ucklichen \emph{Ahnungsverm"ogen}
"uberlassen. Wir sind also leider noch \emph{weit davon entfernt},
da"s wir die Arbeit des Beweisens aufgestellter S"atze etwa einem
\emph{mathematischen Hilfsarbeiter "ubertragen} k"onnen, wie dies
hinsichtlich der Durchf"uhrung einer numerischen Rechnung m"oglich
ist.

Auf den ersten Blick mag freilich die Aufstellung eines allgemeinen,
zwangl"aufigen Verfahrens f"ur das Beweisen mathematischer Aussagen
als ein \emph{hoffnungsloses, ja wahnwitziges Unterfangen}
erscheinen,\,\po{4}\,etwa vergleichbar dem mittelalterlichen Suchen
nach dem \emph{Stein der Weisen}.  Und allerdings w"urde ja die
Erreichung dieses Zieles an \emph{Bedeutung} nicht nur f"ur die
Mathematik, sondern f"ur unsere Erkenntnis "uberhaupt den Wundern, die
man sich von jenem sagenhaften Stein versprach, gewi"s nichts
nachgeben. Es w"are in der Tat, wenn wir genauer zusehen, die ganze
Mathematik in eine \emph{ungeheure Trivialit"at} verwandelt; wir
w"urden jede uns gestellte Aufgabe~--- zwar nicht notwendig in dem
gemeinten Sinne l"osen, was ja unm"oglich sein kann~---, aber doch,
wie Hilbert dies in seinem Pariser Vortrag \del{"uber} "`Mathematische
Probleme"' ausgef"uhrt hat, in dem Sinne erledigen, in welchem es
"uberhaupt einer Erledigung f"ahig ist.

Aber andererseits sieht man \emph{vom rein mathematischen Standpunkt}
in der Tat \emph{gar keine M"oglichkeit}, dieses Problems in seiner
ungeheuren Vielgestaltigkeit irgendwie Herr zu werden, etwa die
\emph{mathematischen Aussagen} je nach dem Grade ihrer Verwicklung
\emph{in Klassen} einzuteilen und auf diese Weise \emph{vom Leichteren
  zum Schwierigeren} fortzuschreiten.

Nichtsdestoweniger haben wir, wie ich behaupten m"ochte, heute bereits
auf einem nicht unmittelbar mathematischen, sondern noch \emph{viel
  grundlegenderen Gebiete} einen derartigen Erkenntnisstandpunkt
erreicht, da"s wir \emph{nicht nur die M"oglichkeit}, sondern auch ein
\emph{gutes Recht} haben, dieses Problem ausdr"ucklich und streng zu
formulieren und in seine Untersuchung einzutreten. In der Tat m"ussen
wir, da es sich hier im Grunde nicht mehr um ein eigentlich
mathematisches, sondern gewisserma"sen um ein
\emph{"`"ubermathematisches"' Problem} handelt, hier \emph{eine Stufe
  tiefer} herabsteigen, und zwar zu der \emph{symbolischen Logik}, wie
sie von \emph{Leibniz} vorausgeahnt, von \emph{Boole} \,\po{5}\, in
Angriff genommen und von \ins{Schr"oder,} \emph{Frege}, \emph{Peano}
und \emph{Russell} m"achtig gef"ordert worden ist.

Tats"achlich gibt sie uns auf ihrem gegenw"artigen
Entwicklungsstandpunkt bereits die M"oglichkeit, das Problem auf eine
\emph{verh"altnism"a"sig einfache und durchsichtige Form} zu bringen,
indem sie uns eine einheitliche Darstellung aller in Betracht
kommenden Aussagen durch ganz wenige Zeichen in die Hand gibt. Wir
k"onnen n"amlich \del{grunds"atzlich} \emph{jede logische oder
  mathematische Aussage} grunds"atzlich \emph{allein mit den folgenden
  wenigen Worten} ausdr"ucken:\footnote{\selectlanguage{english}Ed: The
  following table is originally arranged as follows:
  \foreignlanguage{german}{nicht, oder, alle, (Aussage), Begriff,
    Beziehung, erf"ullt, identisch}. Behmann indicates the corrected
  order by numbers pencilled in as well as a list in the margin.}

\bigskip
\begin{raggedright}
\begin{tabular}{@{}cccccc@{}}
& & & & \multicolumn{2}{c}{\ins{$\overbrace{\hspace{2in}}^{\text{Begriff}}$}}\\
\emph{(Aussage)} & \emph{nicht} & \emph{oder} & \emph{Ding} & \emph{\del{Begriff} \ins{Eigenschaft}} & \emph{Beziehung} \\
$p$ & $\overline{p}$ & $p \, q$ & $a$ & $f^{x}$, $F^{\phi}$ & $f^{xy}$ 
\end{tabular}\hfil\par
\hfill\begin{tabular}{@{}ccc@{}}
\emph{erf"ullt} & \emph{alle} & \emph{identisch} \\
$f_{a}$, $f_{ab}$ & $xf_{x}$ & $a=b$
\end{tabular}
\end{raggedright}
\bigskip

Nat"urlich ist damit nicht gesagt, da"s man in der Praxis tats"achlich
die Zur"uckf"uhrung immer so weit treiben m"usse; vielmehr wird man
auch f"ur andere Begriffe, wie \emph{und}, \emph{es gibt}, \ins{die
  Klasse, die Zahlen und die Verkn"upfung u.s.f.}, soweit es f"ur den
jeweiligen Zweck erw"unscht ist, \emph{neue Zeichen} einf"uhren.

Mit dieser Erkenntnis ist nun immerhin \emph{schon etwas gewonnen}. Wir
brauchen jetzt n"amlich "uberhaupt nur noch \emph{aus den obigen
  Zeichen regelrecht gebildete Ausdr"ucke} in Betracht zu ziehen und
haben uns damit der \emph{Fesseln der Wortsprache} bereits entledigt.

Bekanntlich l"a"st sich die symbolische Logik \emph{axiomatisieren},
d.h. auf ein System verh"alt\-nis\-m"a"sig weniger Grundformeln und
Grundregeln zur"uck\-f"uhren, soda"s auch das \emph{Beweisen von
  S"atzen} nunmehr als ein \emph{blo"ses Rechenverfahren}
erscheint. Man braucht nur noch zu gegebenen Formeln neue
hinzuschreiben, wobei durch Regeln bereits festgelegt ist, was man
jeweils hinschreiben darf. Das Beweisen hat sozusagen den
\emph{Charakter eines Spieles} angenommen. Es ist etwa \,\po{6}\, wie
beim \emph{Schachspiel}, wo man durch \emph{Verschieben} eines der
eigenen Steine, gegebenenfalls mit \emph{Wegnahme} eines gegnerischen,
die \emph{jeweils \del{gegebene} \ins{vorliegende} Stellung in eine
  neue} verwandelt, wobei nur das Verschieben und das Wegnehmen durch
die \emph{Regeln des Spieles} erlaubt sein mu"s.

Aber gerade dieser Vergleich zeigt uns auch in krasser Weise, da"s
uns der eben geschilderte \emph{Standpunkt der symbolischen Logik}
f"ur unser Problem \emph{noch keineswegs gen"ugen} kann. Denn diese
sagt uns wie die Regeln des Schachspiels \emph{nur, was man tun darf},
und \emph{nicht, was man tun soll}. Dies bleibt in dem einen wie in dem
anderen Falle eine Sache des erfinderischen \emph{Nachdenkens}, der
gl"ucklichen \emph{Kombination}. Wir verlangen aber weit mehr: da"s
nicht etwa nur die erlaubten Operationen im einzelnen, sondern auch
der \emph{Gang der Rechnung selbst} durch Regeln festgelegt sein soll,
m.a.W.\ \emph{eine Ausschaltung des Nachdenkens zugunsten des
  mechanischen Rechnens}. Ist irgend eine logische oder mathematische
Aussage vorgelegt, so soll das verlangte Verfahren eine vollst"andige
Anweisung geben, wie man \emph{durch eine ganz zwangl"aufige Rechnung
  nach endlich vielen Schritten} ermitteln kann, ob die gegebene
Aussage \emph{richtig oder falsch} ist. Das eben formulierte Problem
m"ochte ich \emph{das allgemeine Entscheidungsproblem} nennen.

F"ur das Wesen des Problems ist von grunds"atzlicher Bedeutung,
da"s als Hilfsmittel des Beweises \emph{nur das ganz mechanische
  Rechnen} nach einer gegebenen Vorschrift, ohne irgendwelche
Denkt"atigkeit im engeren Sinne, zugelassen wird. Man k"onnte hier,
wenn man will, von \emph{mechanischem} oder \emph{maschinenm"a"sigem
  Denken} reden. (Vielleicht kann man es sp"ater sogar durch eine
Maschine ausf"uhren lassen.)~\po{7}

Es bedarf hier somit nicht blo"s der symbolischen Logik als solcher,
sondern vielmehr einer \emph{besonderen Richtung innerhalb der
  symbolischen Logik}, die nicht nur untersucht, wie man rechnen
\emph{darf}, sondern wie man rechnen \del{mu"s} \ins{\emph{soll}}, um
ein gegebenes Ziel zu erreichen. Diese Richtung m"ochte ich nun aus
einem bestimmten geschichtlichen Grunde als die der \emph{Algebra der
  Logik} bezeichnen, und zwar darum, weil \emph{alle wichtigen
  Vorarbeiten} auf diesem Gebiet bisher gerade unter diesem Namen der
Algebra der Logik (\emph{Boole}, \emph{De Morgan}, \emph{Peirce},
\emph{Schr"oder}) vereinigt sind, \del{und zwar} \emph{vielleicht nur
  zuf"allig} vereinigt sind, denn, soviel ich wenigstens aus
Schr"oders Werk entnehmen konnte, hat wohl keiner der eben Genannten
ernstlich an eine Anwendung auf das Entscheidungsverfahren dabei
gedacht. Aber auch \emph{sachlich} scheint mir dieser Name
\emph{einigerma"sen passend} zu sein, weil dieses engere Gebiet eben
die Ausbildung von \emph{Algorithmen}, d.h. zwangl"aufigen
Rechenverfahren, wie in der Algebra etwa f"ur die Aufl"osung von
Gleichungen, zur Aufgabe haben soll. Ich halte es geradezu f"ur einen
au"serordentlich \emph{gl"ucklichen Umstand}, da"s sich hier eine
Gelegenheit darbietet, die bisher namentlich in der "au"seren Form
ziemlich fragw"urdigen und schwer zug"anglichen, $\langle$um nicht zu
sagen: langweiligen$\rangle$ \ins{eint"onigen} Untersuchungen dieses
Gebietes unter einem \emph{neuen, einheitlichen und wertvollen
  Gesichtspunkt} zusammenzufassen und damit viel tiefe und
beschwerliche Gedankenarbeit vor dem unr"uhmlichen Ende der langsamen
Vermoderung in Bibliotheken zu bewahren. Denn, soweit ich urteilen
kann, hatte sich das \emph{allgemeine Interesse}, selbst bei den
Vertretern der symbolischen Logik, von den Untersuchungen jener alten
Algebra der Logik in den letzten Jahrzehnten so gut wie
\emph{vollst"andig abgewandt.}~\po{8}

Ich brauche wohl nicht erst ausdr"ucklich zu betonen, da"s wir
trotzalledem von der L"osung des vorhin formulierten
Entscheidungsproblems heute noch \emph{au"serordentlich weit entfernt}
sind. Immerhin ist es mir \del{k"urzlich} gelungen, dieses Problem
f"ur einen gewissen \emph{bescheideneren Bereich von
  Aussagen},\footnote{Ed: Noted in left margin:
  \foreignlanguage{german}{\emph{heutiger Stand des Problems}}.} der
  aber nichtsdestoweniger bereits eine recht bedeutende
  Mannigfaltigkeit von Aussagen aufweist, vollst"andig durchzuf"uhren,
  und $\langle$eben dies ist es, wor"uber ich im besonderen zu Ihnen
  sprechen m"ochte.$\rangle$ \ins{gerade dar"uber m"ochte ich in
    diesem Zusammenhang noch etwas sagen.}

\del{Allerdings stehe ich dabei vor einer gewissen \emph{praktischen
    Schwierigkeit.} Es handelt sich n"amlich hier um eine Darstellung,
  die sich auf einem noch \emph{neuen und ungewohnten Gebiet} bewegt
  und mit \emph{v"ollig anderen Mitteln} arbeitet als man sonst in der
  Mathematik gewohnt ist. $\langle$Dazu kommt, da"s die Untersuchung
  ohnehin, zumal sie, was auf anderen Gebieten der Mathematik zumeist
  nicht n"otig ist, \emph{von Grund auf anfangen} mu"s, eine nicht
  unerhebliche Ausdehnung angenommen hat, die sich der verf"ugbaren
  Zeit nur sehr schlecht anpa"st.}
Daher ziehe ich es vor, hier \emph{nur den Gedankengang}
$\langle$meiner Untersuchung$\rangle$ anzugeben, ohne dabei auf
Strenge und Vollst"andigkeit der Darstellung allzuviel Gewicht zu
legen. \del{$\langle$Auch f"ur das Verst"andnis scheint mir dieser Weg
  den Umst"anden nach der f"orderlichste zu sein.$\rangle$}

Zun"achst m"ochte ich einiges "uber die $\langle$von mir$\rangle$
verwendete \emph{Symbolik} sagen. Die Symbolik der \emph{alten Algebra
  der Logik}, selbst diejenige von \emph{Schr"oder}, ist f"ur den
gegenw"artigen Zweck \emph{vollst"andig unbrauchbar}. Es empfiehlt
sich indessen auch nicht, die Zeichengebung der \emph{heutigen
  symbolischen Logik}, etwa diejenige von \emph{Whitehead} und
\emph{Russell},\,\po{9}\, einfach zu "ubernehmen, da diese \emph{wohl
  dem weiteren Zweck} der symbolischen Logik, aber \emph{nicht so sehr
  dem engeren der Algebra der Logik}, wie ich sie verstehe, angepa"st
ist. Daher habe ich mich zwar vorwiegend an die Symbolik von Whitehead
und Russell angeschlossen, aber diese bis zu einem gewissen Grade nach
meinen besonderen Zwecken umgestellt.

Die \emph{Negation} einer Aussage bezeichne ich, wie gesagt, durch
"Uberstreichen: ${\overline {p}}$, die \emph{Disjunktion} als
$p\,q$. Dann stellt sich "`\emph{$p$ und $q$}"' als $\overline
{{\overline{p} \, {\overline{q}}}}$ und "`\emph{Wenn $p$, so $q$}"'
als ${\overline{p}} \, q$ dar. Der Bequemlichkeit halber f"uhre ich
f"ur "`$p$ und $q$"' noch die k"urzere Schreibung $p.q$ ein. Da"s ein
\emph{Ding} (Individuum) $a$ \emph{unter einen Begriff $f^{x}$
  f"allt}, schreibe ich so: $f_a$. Ist die \emph{unbestimmte Aussage
  $f_x$ f"ur jedes Ding $x$ richtig}, so schreibe ich, um dies
anzudeuten: $x\,f_x$. Da"s \emph{mindestens ein Ding} die Aussage
$f_x$ erf"ullt, l"a"st sich dann als
$\overline{x\,{\overline{f_{x}}}}$ ausdr"ucken. Ich schreibe daf"ur
"uberdies, die "ubereinanderstehenden Teile der beiden Striche
gewisserma"sen gegeneinander weghebend: ${\overline{x}\,f_{x}}$.

\emph{Klammern} benutze ich im wesentlichen nur, um die Verkn"upfungen
der Disjunktion und der Konjunktion voneinander zu trennen (falls die
letzte durch den Punkt bezeichnet wird), sonst nur \emph{zur
  Hervorhebung}.

So z.B. $x\,f_{x}g_{x}$\ , $x\,f_{x}p$. Das letzte ist \ins{an
sich} \emph{zweideutig!} Es gilt aber $(xf_{x})p
\leftrightarrow x(f_{x}p)$ Regel II$^{\ast}$. Entsprechend schreibe
ich $x.f_{x}.g_{x}$ \ \ ${\overline{x}}.f_{x}.p$, wo der erste Punkt
nur Lesezeichen ist.

Unter einen Begriff k"onnen nun \emph{selbst wieder Begriffe}
fallen. Ist $F^{\phi}$ ein \emph{Begriff zweiter Stufe}, so schreibe
ich $F_{\phi}$ f"ur die Aussage, da"s der Begriff $\phi$ unter den
Begriff $F$ f"allt. Die Bedeutungen von $\phi F_{\phi}$ und
${\overline{\phi}} F_{\phi}$ sind hiernach ohne weiteres klar.  

Ich komme nun zur \emph{Erkl"arung der Rechenregeln}, die \po{10}
dar"uber Auskunft geben, welche Operationen man mit den
hingeschriebenen Zeichen vornehmen darf. Solche Regeln sind
\emph{auch in der axiomatisch betriebenen symbolischen Logik}
bekanntlich unentbehrlich. W"ahrend man aber dort ihre Zahl und den
Gehalt jeder einzelnen auf ein Minimum zu beschr"anken sucht,
verfolgen wir hier \emph{das gegenteilige Ziel}, nicht: ihre Zahl
beliebig zu vermehren~--- das w"urde das Ged"achtnis unn"otig
belasten~---, wohl aber, jeder einzelnen als auch dem Regelsystem
insgesamt eine \emph{m"oglichst gro"se praktische Tragweite} zu
geben. Im Gegensatz zu[r] axiomatischen Untersuchung sind hier also
\emph{durchaus praktische R"ucksichten} ma"sgebend.
\bigskip

\begin{enumerate}
\item[I.] Satz von der doppelten Verneinung. \qquad
  $\overline{\overline {p \, {\overline q}}} \, r \leftrightarrow (p
  \, {\overline q})r$

\item[II, II$^{\ast}$.] Vereinigungssatz. \qquad $(p \, {\overline
  q})r \leftrightarrow p ({\overline q} \, r) \leftrightarrow p \,
  {\overline q}\, r$, \ \ $(xf_{x})p \leftrightarrow x(f_{x}p)$

\item[III, III$^{\ast}$.] Vertauschungssatz.  
\[
\begin{array}[t]{@{}r@{\,}lr@{\,}l@{}}
{\overline{p \, q}} \, r \, {\overline{s}} & \leftrightarrow {\overline{s}} \,
  {\overline{q \, p}} \, r & x {\overline y}f_{x}g_{y} &
  \leftrightarrow xf_{x}{\overline y}g_{y} \\
xy\,f_{xy} & \leftrightarrow yx\,f_{xy}, & {\overline x}\,{\overline y}f_{xy} & \leftrightarrow {\overline y}\,{\overline x}f_{xy}
\end{array}
\]

\item[IV, IV$^{\ast}$.] Verschmelzungssatz. 
\[
p \, p \, q \leftrightarrow p\,q, \qquad
\left\{\begin{array}{@{}r@{\,}l@{}}
x\,f_{x}.x\,g_{x} & \leftrightarrow x.f_{x}.g{x} \\
xf_{x}\,xg_{x} & \leftrightarrow xf_{x}g_{x}
\end{array}\right.
\]

\item[V, V$^{\ast}$.] Eliminationssatz.  \qquad [$p \, {\overline{{\overline{q} \, r}}}.{\overline{q}}\, r \, {\overline{p \, s}} \rightarrow p \, {\overline{p \, s}}$] 

\item[V, V$^{\ast}$.] Implikationssatz.

\item[VI, VI$^{\ast}$.] Satz von der Ausf"uhrung der Negation.
\[
\left\{\begin{array}{@{}r@{\,}l@{}}
{\overline{p \, q}} & \leftrightarrow {\overline{p}}.{\overline{q}},  \\
{\overline{xf_{x}}} & \leftrightarrow {\overline x} {\overline{f_{x}}} 
\end{array}
\right.
\qquad
\left\{\begin{array}{@{}r@{\,}l@{}}
{\overline{p.q}} & \leftrightarrow {\overline{p}} \, {\overline{q}} \\
{\overline{{\overline x}f_{x}}} & \leftrightarrow x {\overline{f_{x}}}
\end{array}\right.
\]

\item[VII.] Distributionssatz. \qquad $p \, q.r \leftrightarrow (p.r)(q.r)$, \ \ $(p.q)r \leftrightarrow p \, r.q \, r$

\item[VIII.] Dualit"atssatz.

\item[IX.] Vereinfachungss"atze.

\item[X.] Verdr"angungssatz.

\item[XI.] Satz von der Umschreibung der singul"aren Aussage. \\
\hspace*{2em} $f_{a} \leftrightarrow x \, {\overline{x = a}}f_{x}
\leftrightarrow {\overline x}. x=a.f_{x}$ \po{11}
\end{enumerate}

Diese wenigen Andeutungen m"ogen gen"ugen.

Bez"uglich unseres allgemeinen Zieles erinnere ich daran, da"s f"ur
einen gewissen ganz einfachen Bereich von Aussagen, n"amlich alle
diejenigen, die \emph{nur die Worte nicht und oder} ("uberdies, wenn
man will, \emph{und, folgt, "aquivalent}) voraussetzen, die
\emph{L"osung des Entscheidungsproblems} bereits seit einiger Zeit
\emph{bekannt} ist.\footnote{Ed.: The connective words in this
  sentence are individually underlined in black ink.} Sogar
\emph{zwei} L"osungen.

Ist z.B. die Aussage $p \supset \colon p \supset q
\mathrel{.{\supset}.} q$ vorgelegt, die ich so schreibe: ${\overline
  p} \, {\overline{{\overline p} \, q}} \, q$, so braucht man hier
f"ur $p$ und $q$ nur die Werte \del{"`richtig"'} \ins{"`\emph{wahr}"'}
($\curlyvee$) und "`falsch"' ($\curlywedge$) auf die \emph{vier
  m"ogliche\ins{n} Arten} einzusetzen, wobei
\[
\overline{\curlyvee} \leftrightarrow \curlywedge,\qquad
\overline{\curlywedge} \leftrightarrow \curlyvee,\qquad
[\curlyvee\curlyvee \leftrightarrow \curlyvee\curlywedge \leftrightarrow \curlywedge\curlyvee \leftrightarrow \curlyvee],\qquad
\curlywedge\curlywedge \leftrightarrow \curlywedge
\]
und erh"alt jedesmal den Wahrheitswert $\curlyvee$, womit die
Richtigkeit, genauer: die \emph{Allgemeing"ultigkeit} der gegebenen
Aussage best"atigt ist.

Eleganter ist das \emph{Verfahren der Normalform.} Verm"oge des
Distributionssatzes u.s.f. kann man jede Aussage der betrachteten Art
als \emph{Konjunktion von lauter Disjunktionen von Bestandteilen $p$,
  ${\overline p}$} schreiben.
\[
{\overline p}(p.{\overline q}) q \leftrightarrow {\overline p} \, p
\, q.p \, {\overline q} \, q
\]
Hier enth"alt das erste Glied den Bestandteil ${\overline p} \, p
\leftrightarrow \curlyvee$ und das zweite den Bestandteil ${\overline
  q} \, q \leftrightarrow \curlyvee$. Diese Bedingung ist andererseits
\emph{auch hinreichend}, soda"s man jede derartige Aussage auf diese
Weise pr"ufen kann.

Ich hatte anf"anglich versucht, das \emph{erste Verfahren auf h"ohere
  Aussagen zu "ubertragen.} Man kann n"amlich auf Grund unserer
Regeln \emph{alle Operatoren nach vorn} bringen, soda"s rechts eine
Verkn"upfungsaussage allein aus \emph{nicht, oder, und}
"ubrigbleibt. Z.B.
\[
x{\overline{y{\overline{f_{xy}}}g_{y}.{\overline z}{\overline{h_{yz}}}}} \ \leftrightarrow \ x{\overline y}z {\overline{\overline{f_{xy}}g_{y}.{\overline{h_{yz}}}}} \ \leftrightarrow \ x{\overline y}z\,\textrm{(Normalform)} 
\]
\hbox to\parindent{\pind}\makepagination{12}Dies f"uhrt jedoch, soviel
ich sehe, \emph{nicht zum Ziel}. Daher wird man versuchen, das
\emph{Verfahren der Normalform} zu verallgemeinern.

Zun"achst m"ochte ich den \emph{Bereich von Aussagen} genauer
umgrenzen, f"ur den ich das Entscheidungsproblem durchgef"uhrt
habe. Was diejenigen Aussagen betrifft, die "uberhaupt f"ur unser
Problem in Frage kommen, so werden wir sagen, da"s
\emph{au"serlogische konstante Bestandteile}, wie $a$, $f^{x}$,
$f^{xy}$ gewi"s \emph{nicht} vorkommen werden. Nur
\emph{ver"anderliche}, d.h. gleichzeitig durch Operatoren vertretene:
\[
x,\quad \phi^{x},\quad \phi^{xy},\quad \Phi^{\phi},\text{
  u.s.f., au"serdem aber } x=y.
\]
Wir beschr"anken uns zun"achst auf: $x$, $\phi^{x}$, lassen also nur
die folgenden Worte zu: \emph{nicht, oder, alle, Ding,}
\del{\emph{Begriff}} \ins{\emph{Eigenschaft}} (von Dingen),
\emph{erf"ullt}. Nat"urlich mu"s jedes Dingsymbol und jedes
Begriffssymbol irgendwo zu Anfang oder innerhalb der Aussage
\emph{zugleich als Operator} $x$, ${\overline x}$, $\phi$, ${\overline
  \phi}$ vorkommen. \ins{Auch von $x=y$ sehen wir vorl"aufig ab.}
\emph{Teilbereich der Aussagen zweiter Ordnung} nach Russell.

Wir betrachten das triviale Beispiel: $\phi{\overline \chi} x
{\overline {\phi_{x}}} {\overline {\chi_{x}}}$, etwa zu lesen: "`Zu
jedem Begriff $\phi$ gibt es einen Begriff $\chi$, soda"s, was
unter $\phi$ f"allt, nicht unter $\chi$ f"allt."' Oder:
\[
\phi \chi \psi \ {\overline{x\, {\overline {\phi_{x}}} \chi_{x}}} \  
{\overline{x\, {\overline {\chi_{x}}} \psi_{x}}} \  
x\, {\overline {\phi_{x}}} \psi_{x}
\]
Barbara. \emph{Besonderer Fall}. Die Begriffsoperatoren sind s"amtlich
allgemein und am Anfang der Aussage. Sie k"onnen auch irgendwie
innerhalb der Aussage verstreut sein \ins{und
partikul"ar!} Diesen Aussagenbereich nenne ich den den
\emph{Bereich B}.

\emph{Allgemeines Programm.} Die geg. Aussage ist nach einer gewissen
Vorschrift in eine solche umzuschreiben, die \emph{einen Begriff
  weniger} enth"alt. Dies wird so lange fortgesetzt bis \emph{alle
  Begriffe $\phi$, $\chi$, \dots\ verschwunden sind} und eine ganz
triviale Aussage "ubrig bleibt, f"ur welche die Entscheidung ohne
weiteres vollzogen werden kann.~\po{13}

\emph{Eliminationsverfahren.} (Aussonderung ?) \emph{Formulierung.}
Ich nehme den \emph{am weitesten rechts stehenden Begriffsoperator und
  den zugeh"origen Operand:}
\[
\phi F_{\phi \chi u v} \qquad\text{oder}\qquad {\overline \phi} F_{\phi \chi u v}.
\]
$\chi$ hat den Operator weiter links, ebenso $u$ und $v$, au"serdem
k"onnen aber Operatoren $x$, ${\overline y}$ innerhalb $F$
auftreten. Diese lie"sen sich aber oben nicht mit
auff"uhren. Dagegen enth"alt $F$ nach Voraussetzung \emph{keine
  Begriffsoperatoren}.

Wir sehen nun auch noch von den Operatoren $\phi$ bez. ${\overline
  \phi}$ ab und beschr"anken uns \emph{allein auf den Operand}
$F_{\phi \chi \ldots u v \ldots}$, um diesen weiter umzuformen.

\emph{Wesentlicher Gedanke:} Bei der Betrachtung von $F$ braucht man
sich um die \emph{au"serhalb} stehenden Operatoren \emph{nicht} zu
k"ummern. $\phi$, $\chi$ sind daher als \emph{konstante Begriffe},
$u$, $v$ als \emph{konstante Dinge} zu behandeln.

Wir haben also die Aufgabe vor uns, den \emph{Begriff $\phi$ zu
  eliminieren} aus einer Aussage ${\overline \phi} \, Q$. (Es
gen"ugt, uns auf diesen Fall zu beschr"anken, da der andere durch
Verneinung entsteht.) Hier ist nun $Q$ eine Aussage, die
m"oglicherweise \emph{konstante und ver"anderliche Dinge}, aber
keine anderen als \emph{konstante Begriffe}, d.h. eben keine
Begriffsoperatoren, enth"alt. Alle derartigen Aussagen fasse ich
zusammen zu dem \emph{Aussagenbereich $A$}. Teilbereich der
\emph{Aussagen erster Ordnung}.

Und diese Aussagen des Bereiches $A$ haben nun die angenehme
Eigenschaft, da"s sie sich \emph{tats"achlich als Normalformen}
schreiben lassen. Die Frage ist nun: Aus welchen \emph{Bausteinen}
werden sich diese Normalformen zusammensetzen. (Fr"uher: $p$,
${\overline p}$) Es sind zun"achst, wenn ich die Konstanten
\del{wieder durch lateinische Buchstaben} \ins{als
solche} bezeichne, solche von der Form $f_a$,
weiterhin:
\[
\left\{
\begin{array}{lll}
xf_{x}, & x\,f_{x}g_{x}, & x\,f_{x}g_{x}h_{x}, \dots\\ {\overline
  x}.f_{x}, & {\overline x}.f_{x}.g_{x}, & {\overline
  x}.f_{x}.g_{x}.h_{x}, \dots
\end{array}\right.
\]
\hbox to\parindent{\pind}\makepagination{14}(Die $f$, $g$, $h$ k"onnen
auch \emph{zum Teil verneint} sein.) Andere Formen
kommen nicht vor. In der Tat ist z.B. nach IV$^{\ast}$: 
\[
{\overline x} f_{x} g_{x} \ \ \leftrightarrow \ \ {\overline x} f_{x}
{\overline x} g_{x}
\]

F"ur die obigen Bausteine f"uhre ich eine \emph{abgek"urzte Schreibung} ein:

\emph{Neue Klassensymbolik}. \qquad $\alpha$, $\dot{\alpha}$, $\vee$,
$\wedge$, $\alpha_{a}$.

\def\overbow#1#2{%
\tikz[baseline=(X.base)]{\node[inner xsep=0pt,inner ysep=1pt] (X) {$#1$};
\foreach \x in {1,...,#2}
\draw (current bounding box.north west) .. controls +(30:.12) and +(150:.12) .. (current bounding box.north east);}}
\def\underbow#1#2{%
\tikz[baseline=(X.base)]{\node[inner xsep=0pt,inner ysep=1pt] (X) {$#1$};
\foreach \x in {1,...,#2}
\draw (current bounding box.south west) .. controls +(-30:.12) and +(-150:.12) .. (current bounding box.south east);}}

Bausteine
\[
\raisebox{-2ex}{$\biggl\{$}
\begin{array}[t]{@{}llll}
\underbow{\alpha}{1} & \underbow{\alpha \beta}{1} & \underbow{\alpha \beta
  \gamma}{1} &\ldots\\[1ex] \overbow{\alpha}{1} & \overbow{\alpha
  \beta}{1} & \overbow{\alpha \beta \gamma}{1} &\ldots
\end{array}
\]
\emph{Bedeutungen:} \qquad Beispiel: 
\begin{align*}
x {\overline{y}}.f_{x}g_{y}.h_{x}k_{y} & \leftrightarrow\\
& \leftrightarrow  (xf_{x})({\overline{y}}g_{y}).(xf_{x}h_{x})({\overline{y}}.g_{y}.k_{y}).(xh_{x})(yk_{y}) \\
& \qquad\qquad\text{\quad durch zwangsl"aufige Umformung}\\
& \leftrightarrow \underbow{\alpha}{1}\, \overbow{\beta}{1} . \underbow{\alpha\gamma}{1}\,\overbow{\beta\delta}{1} . \underbow{\gamma}{1}\,\overbow{\delta}{1}.
\end{align*}

Unser Problem lautet nunmehr ${\overline{\rho}}\,F_{\rho\alpha\beta
  \ldots ab \ldots}$. Wir haben also \emph{links einen Operator
  $\overline{\rho}$} und rechts davon einen Ausdruck, nehmen wir an,
eine \emph{disjunktive Normalform}, d.h. eine Disjunktion aus lauter
Konjunktionen der obigen Bestandteile. Da wir nach IV$^\ast$
\emph{$\overline{\rho}$ zu jedem einzelnen Disjunktionsglied} setzen
k"onnen, haben wir also nur noch Bestandteile~$\overline{\rho}\,
\mathfrak{S}$ zu betrachten, wo $\mathfrak{S}$ eine \emph{Konjunktion
  von Bausteinen} ist. Und zwar gen"ugt es hier im besonderen, die
Form 
\[
\overline{\rho}.\underbow{\alpha\rho}{1}.\underbow{\beta\dot\rho}{1}\, .\,
\overbow{\gamma\rho}{1}\, .\, \overbow{\gamma'\rho}{1}\ldots
\overbow{\delta\dot{\rho}}{1}.\overbow{\delta'\dot{\rho}}{1}\ldots
\]
zu betrachten, wie in der Hauptsache bereits der \emph{alten Algebra
  der Logik} bekannt war. Freilich war diese nicht von dem
Bereich~$A$, sondern von dem sehr viel engeren Bereich der \del{durch
  ihre Symbolik, d.h.}  der \emph{durch blo"se Klassenverkn"upfungen}
herstellbaren Aussagen ausgegangen, wie sie gar nicht anders konnte,
da ihre Symbolik die Darstellung der Aussagen unseres Bereiches~$A$
(mit beliebig durcheinander gew"urfelten Individuenoperatoren)
"uberhaupt nicht zulie"s. Des obige Problem hat \emph{Schr"oder zu
  l"osen vesucht}; es ist ihm aber trotz aller\, \po{15}\,M"uhe
\emph{nicht gelungen}. Er ist "uber gewisse spezielle Ergebnisse nicht
hinausgekommen. Der Grund liegt nach meiner Ansicht durchaus in der
\emph{Spr"odigkeit seiner Symbolik}.

Ein \emph{besonders einfacher Fall} des obigen Problems ist der folgende:
\[
{\overline \rho}\, .\, \underbow{\alpha\rho}{1}\, .\, \underbow{\beta\dot{\rho}}{1}  \leftrightarrow  \underbow{\dot\alpha\beta}{1}
\]
Oder bei anderer Wahl der Gr"o"sen: 
\[
{\overline \rho}\, .\,\underbow{\dot{\alpha}\rho}{1}\, .\,
\underbow{\dot{\rho}\beta}{1} \leftrightarrow \underbow{\dot{\alpha}
  \beta}{1}
\]

Dies ist nichts anderes als der bekannte Schlu"s
\emph{Barbara}. $\underbow{\dot{\alpha} \rho}{1}$ hei"st n"amlich,
da"s $\alpha$ eine Teilklasse von $\rho$ ist. Gibt es zu $\alpha$ und
$\beta$ eine "`Zwischenklasse"' $\rho$, so ist gewi"s auch $\alpha$
eine Teilklasse von $\beta$.  Aber auch das Umgekehrte gilt. Ist
$\alpha$ eine Teilklasse von $\beta$, so l"a"st sich gewi"s auch eine
Zwischenklasse angeben. Denn es ist ja $\alpha \subset \alpha \subset
\beta$ und $\alpha \subset \beta \subset \beta$.

Dies ist das \emph{grundlegende Eliminationsergebnis}, auf das alle
Eliminationen letztlich zur"uckgef"uhrt werden. Tats"achlich kommt die
Klasse $\rho$ bez. der Begriff $\phi$ \emph{rechts nicht mehr vor.}

Von hier aus gelingt es nun durch einen \emph{einfachen Gedanken} das
obige allgemeine Problem Schr"oders \emph{rein rechnerisch zu
  erledigen}, allerdings durch einen Gedanken, den Schr"oder selbst
unm"oglich haben konnte, weil er ganz wesentlich auf der Verwendung
der \ins{\emph{modernen}} \emph{verbesserten
  Symbolik} beruht. Es gilt n"amlich, (als Sonderfall von
III$^{\ast}$) die logische Regel:
\[
{\overline x} \, {\overline y}\,f_{xy} \leftrightarrow  {\overline y} \, {\overline x}\, f_{xy}.
\]
Ebenso, wenn Begriffsoperatoren vorkommen: 
\[
{\overline \phi}\,
{\overline x}\,F_{\phi x} \leftrightarrow {\overline x}\,
{\overline \phi}\, F_{\phi x}.
\]

Nun stecken aber in den $\gamma$- und den $\delta$-Gliedern
\emph{partikul"are Operatoren} ${\overline x}$, ${\overline y}$, und
diese kann man verm"oge III$^{\ast}$ "uber den Klassenoperator
${\overline \rho}$ \emph{nach links hinausschieben}, soda"s wir die
rechts verbleibenden Individuen $x$, $y$ einfach als \emph{konstant}
anzusehen haben, da wir ja nun den Ausdruck von \po{16}
${\overline{\rho}}$ an betrachten werden. Verm"oge XI l"a"st sich
der neue Operand von ${\overline{\rho}}$ \emph{ganz aus allgemeinen
  Aussagen} zusammensetzen und "uberdies wieder die Anzahl der
$\alpha$- und der $\beta$-Glieder \emph{je auf Eins} zur"uckf"uhren,
so da"s die Elimination nach der vorhin angegebenen Regel ohne
weiteres vollzogen werden kann. Das \emph{Ergebnis} ist folgendes:
\[
\alpha \beta\, . \bigl[ \overbow{\gamma \beta}{1}\, .\,
  \overbow{\gamma ' \beta}{1} \ldots \mid \ . \, \overbow{\delta
    \alpha}{1}\, . \, \overbow{\delta ' \alpha}{1} \ldots \bigr]
\text{\quad \emph{Rohe Resultante!}}
\]
Es besage: $\alpha \cup \beta = \vee$, $\overbow{\gamma \beta}{1}$,
$\overbow{\gamma '\beta}{1}$, u.s.w. $\overbow{\delta \alpha}{1}$,
$\overbow{\delta '\alpha}{1}$, u.s.w. Die Klammern und der Strich
sollen die \emph{Zusatzbedingung} zum Ausdruck bringen, da"s die
Elemente, die $\gamma$, $\gamma'$, $\ldots$ mit $\beta$ einerseits und
$\delta$, $\delta'$, $\ldots$ mit $\alpha$ andererseits gemeinsam
haben, \emph{durchweg verschieden w"ahlbar} sein sollen. In
Begriffsschrift:
\[
x \alpha_{x} \beta_{x}\, .\, {\overline u} {\overline{u'}}\_\_
\ {\overline v} {\overline {v'}}\_\_ \, . \, (u u'\_\_ \ \mid \ v
v'\_\_ \ )\, .\, \gamma_{u}\, .\, \beta_{u}\, .\, \gamma'_{u'}\, .\,
\beta_{u'} \ldots \delta_{v}\, .\, \alpha_{v}\, .\, \delta'_{v'}\, .\,
\alpha_{v'} \ldots
\]  
wo $(u u'\_\_ \ \mid \ v v'\_\_ \ )$ nur besagen soll, da"s jedes
$u$ von jedem $v$ verschieden ist.

Damit ist die \emph{erste Begriffselimination nunmehr
  geleistet}. Wollen wir das Eliminationsgesch"aft \emph{fortsetzen},
so hat die Sache allerdings einen \emph{Haken}. Denn die gewonnene
Resultante geh"ort ja ihrerseits \emph{dem Bereich $A$ im allgemeinen
  nicht mehr an}, da sie ja noch die \emph{Beziehung der Identit"at}
enth"alt, l"a"st sich somit auch nicht als Verkn"upfungsaussage aus
den vorhin aufgez"ahlten Bausteinen und daher auch nicht als
Normalform aus solchen schreiben.\footnote{Ed: The left margin
  contains the text \foreignlanguage{german}{"`\emph{Problem der
      Klausel!}"'} with an arrow pointing to the beginning of this
  paragraph.}

Es empfiehlt sich nun, um die dem Problem angemessene Allgemeinheit zu
erreichen, die \emph{Identit"at} von vornherein \emph{mit unter die
  Bestandteile der Aussagen von $A$} aufzunehmen. Den so erweiterten
Bereich nenne ich den \emph{Bereich} $A^{\ast}$.

Dies Problem erschien mir \emph{anfangs hoffnungslos verwickelt.} Nach
manchem Herumprobieren gelang es mir indessen, die \po{17} Aussagen
von $A^{\ast}$ als \emph{Normalformen} darstellbar zu erweisen aus den
\emph{fr"uheren Bestandteilen} und solchen von den Formen:
\begin{align*}
\underbow{\alpha}{2} \qquad \underbow{\alpha\beta}{2} \qquad \underbow{\alpha\beta\gamma}{2} \qquad \ldots \qquad \underbow{\alpha}{3} \qquad  \underbow{\alpha\beta}{3} \qquad \ldots \\
\overbow{\alpha}{2} \qquad \overbow{\alpha\beta}{2} \qquad \overbow{\alpha\beta\gamma}{2} \qquad \ldots \qquad \overbow{\alpha}{3} \qquad  \overbow{\alpha\beta}{3} \qquad \ldots
\end{align*}

Hier bedeutet $\overbow{\alpha}{2}$, da"s $\alpha$ mindestens 2
Individuen enth"alt $\longrightarrow$.
Z.B. $\underbow{\dot{\alpha}\beta}{2}$, da"s $\alpha$ alle Elemente
von $\beta$ mit h"ochstens einer Ausnahme enth"alt.

Die \emph{Elimination} aus einem solchen Ausdruck geht nun
\emph{bemerkenswerterweise genauso glatt} vonstatten wie im Falle des
Bereichs~$A$, soda"s wir nunmehr aus einer Aussage des
\emph{erweiterten Bereichs}~$B^{\ast}$ tats"achlich alle
Begriffsoperatoren der Reihe nach eliminieren k"onnen. Was nach allen
Zur"uckf"uhrungen \emph{schlie"slich "ubrig} bleibt, kann
z.B. ein Ausdruck wie $\underbow{\vee}{1}$ sein, der die gegebene Aussage als
\emph{unbedingt}, d.h. in jedem Individuenbereich \emph{g"ultig}
erweist~--- solche werden uns nat"urlich vor allem interessieren, oder
aber etwa $\overbow{\wedge}{1}$ was eine \emph{unbedingt falsche} Aussage
bedeutet. Im allgemeinen Fall wird dagegen ein \emph{Ausdruck nach Art
  des folgenden} herausspringen, (den ich absichlich m"oglichst
verzwickt aufgebaut habe):
\[
\overbow{\vee}{5}\,(\overbow{\vee}{2}\,.\,\underbow{\wedge}{4})(\overbow{\vee}{3}\,.\,\underbow{\wedge}{2})\,\underbow{\wedge}{1}
\]
Dies ist eine \emph{disjunktive Normalform}. 

\emph{Bedeutung der Glieder!} 

Die gegebenen Aussage ist im allgemeinen richtig, aber falsch, wenn es
genau 1 Individuum oder wenn es genau 4 Individuen gibt.

\ins{\emph{Bedeutung z.B. f"ur die Geometrie!}}\footnote{Ed: A paragraph
  in shorthand omitted.}
\end{document}